\documentclass{article}

\newcommand{\nc}{\newcommand}\nc{\á}{\´{a}} \nc{\â}{\^a}  \nc{\ã}{\~a} \nc{\é}{\'e} \nc{\ê}{\^e}\nc{\ia}{\'{\i}}
\nc{\í}{\ia} \nc{\ó}{\'o} \nc{\ô}{\^o} \nc{\õ}{\~o}\nc{\ú}{\'{u}}
\nc{\al}{\mbox{$\alpha$}} \nc{\ap}{\approx} \nc{\bc}{\begin{center}}
\nc{\ec}{\end{center}} \nc{\beq}{\begin{equation}}
\nc{\eeq}{\end{equation}} \nc{\beqn}{\begin{eqnarray}}
\nc{\eeqn}{\end{eqnarray}} \nc{\beqne}{\begin{eqnarray*}}
\nc{\eeqne}{\end{eqnarray*}} \nc{\be}{\mbox{$\beta$}}
\nc{\bp}{\bar{\mbox{$\phi$}}} \nc{\ca}{\c c\~ao\ } \nc{\cs}{\c c\~oes\ } 
\nc{\cl}{\mbox{$\cal L$}} \nc{\ch}{\mbox{$\cal H$}}
\nc{\chr}{\mbox{$\cal R$}} \nc{\cf}{\mbox{$\cal F$}}
\nc{\cu}{\mbox{$\cal U$}} \nc{\da}{\dot{A}} \nc{\dq}{\dot{q}}
\nc{\op}{\dot{p}} \nc{\ddq}{\ddot{q}} \nc{\dr}{\dot{\!\vr}}
\nc{\pp}{\dot{\!\vp}} \nc{\de}{\delta} \nc{\eq}{equation\ }
\nc{\Eq}{Eq.\ } \nc{\eqs}{equations\ } \nc{\Eqs}{Eqs.\ }
\nc{\edp}{partial differential equation\ } \nc{\edps}{partial
differential equations\ } \nc{\f}{function\ } \nc{\fs}{functions\ }
\nc{\fa}{\frac{1}{2}} \nc{\fr}{\frac{d}{dt}}
\nc{\fx}[1]{\frac{d}{d\:#1}} \nc{\fh}{Hamiltonian function\ }
\nc{\fl}{Lagrangian function\ } \nc{\fle}{\hs{.5cm}\Longrightarrow
\hs{.5cm}} \nc{\ga}{\mbox{$\gamma$}} \nc{\Ga}{\mbox{$\Gamma$}}
\nc{\h}{Hamiltonian\ } \nc{\hr}{Hamiltonizate\ }
\nc{\hj}{Hamilton-Jacobi\ } \nc{\hz}{Hamiltonization\ }
\nc{\hf}{\hfill} \nc{\hp}{\hs\parindent} \nc{\hs}{\hspace}
\nc{\la}{Lagrangian\ } \nc{\e}{\left}\nc{\ri}{\right}\nc{\na}{\vec{\nabla}} \nc{\no}{\nonumber}
\nc{\noi}{\noindent} \nc{\p}{\mbox{$\phi$}} \nc{\pr}[1]{\partial#1}
\nc{\paf}[3]{\frac{\pr#1}{\pr#2_#3}} \nc{\pa}[3]{{\pr#1}/{\pr#2_#3}}
\nc{\pac}[3]{\frac{\pr#1}{\pr#2^#3}}\nc{\pal}[3]{\frac{\pr^2#1}{\pr#2\pr#3}}
\nc{\pu}[4]{\frac{\pr#1^#2}{\pr#3^#4}}
\nc{\pub}[4]{\frac{\pr#1_#2}{\pr#3_#4}}
\nc{\puf}[4]{{\pr#1_#2}/{\pr#3_#4}}
\nc{\po}[4]{\frac{\pr#1_#2}{\pr#3^#4}}
\nc{\ps}[5]{\frac{\pr^2#1}{{\pr#2_#3}{\pr#4_#5}}}
\nc{\psd}[5]{{\pr^2#1}/{{\pr#2_#3}{\pr#4_#5}}}
\nc{\pe}[2]{\frac{\pr#1}{\pr#2}} \nc{\pee}[2]{\pr#1/\pr#2}\nc{\ped}[2]{\frac{\pr^2#1}{\pr#2^2}}\nc{\pro}{procedure\ }\nc{\q}{\quad}\nc{\qq}{\qquad}
\nc{\ql}{\linebreak} \nc{\qf}{\pagebreak} \nc{\s}{s\~ao\ }
\nc{\Th}{\mbox{$\Theta$}} \nc{\vr}{\vec{\:r}} \nc{\vp}{\vec{\:p}}\nc{\vx}{\vec{\:X}}
\nc{\vs}{\vspace} \nc{\va}[3]{{\de#1}/{\de#2_#3}}
\nc{\vac}[3]{\frac{\de#1}{\de#2^#3}}
\nc{\vaf}[3]{\frac{\de#1}{\de#2_#3}}
\nc{\vef}[2]{\frac{\de#1}{\de#2}} \nc{\ve}[2]{{\de#1}/{\de#2}}
\begin{document}

\title{\bf GENERAL SOLUTION TO UNIDIMENSIONAL HAMILTON-JACOBI EQUATION}

\author{{\bf Maria Lewtchuk Espindola}\\Universidade Federal da Para\ia ba,
\ DM/CCEN\\58051-970, \ Jo\ão Pessoa, \ PB, \
Brazil\thanks{mariia@mat.ufpb.br}}
\date{02/10/2010}
 \thispagestyle{empty} \maketitle\vs{.7cm}
\begin{abstract}
A method for finding the general solution to the partial differential equations: \ $F(u_x,u_y)=0$; \ $F(f(x)\:u_x,u_y)=0$ \ (or \ $F(u_x,h(y)\:u_y)=0$) \ is presented, founded on a Legendre like transformation and a theorem for Pfaffian differential forms. As the solution obtained depends on an arbitrary function, then it is a general solution. As an extension of the method it is obtained a general solution to PDE: \ $F(f(x)\:u_x,u_y)=G(x)$, and then applied to unidimensional Hamilton-Jacobi equation.\end{abstract}\vs{.7cm}


\hs{2.4cm}Nonlinear PDEs;

\hs{2.4cm}Hamilton-Jacobi Equation.\vs{.2cm}

\section{Introduction}

\hp The practical and conceptual importance of the Hamilton-Jacobi
\eq can be extensively pointed: as a fundamental concept in
classical mechanics \cite{ARN}; as a practical tool for solving
differential \eq \cite{CHO}; as a base to quantization \cite{DU};
as an approximation of zero order in the WKB method \cite{SCHI}; ...

The solutions of the \hj  \eq are usually determined as integral solutions
through the method of separation of variables. But general solutions of this \eq are more important either
by its conceptual  meaning \cite{DIR}, as by the characteristic of a general solution (an infinity of integral solutions).

Unfortunately, due to the nonlinearity of this \eq till now there is no available technique to determine a general solution \cite{SNE,FOR} in most problems.

The purpose of this article is solve this centenary problem. The method solves the unidimensional problem changing an insurmountable problem by an eventual practical problem of solving an algebraic equation. Although to this algebraic problem the actual advanced computational numerical techniques can be applied.

The procedure applied to solve this problem ia an extension of that developed for PDEs, linear or not, of one of the types: \[F(p,q)=0; \q F(f(x)\:p,q)=0 \q (or \q F(p,h(y)\:q)=0),\] \noi where \
$p=\pee{u}{x}$, \ $q=\pee{u}{y}$ \ and \
$u=u(x,y)$ \cite{ENA,CNM}.

First this method of obtaining a general solution will be summarized, then its extension to the \hj equation will be  developed.

\section{General Solution to the PDE \ $\mathbf{F(p,q)=0}$}

\hp There are few books that accost the subject of partial differential equations of first order and also the methods presented as the Charpit's, or the characteristic's, or separation of variables one, or another techniques only supplies  integral (complete) solutions to this type of equation unless it is a linear one (Lagrange PDE) \cite{SNE,FOR,COUR,IOR}. \ In the method developed an infinity of integral solutions are obtained. Because it furnishes a general solution in an implicit form and in some cases a explicit one. In the non linear PDEs the general solution is really obtained exchanging an impossible problem by the solution of algebraic equations.

The procedure developed is based on a Legendre like transformation and the use of a theorem that gives the condition of integrability for Pfaffian differential forms \cite{SNE}.

Consider a PDE of first order \beq F(p,q)=0,\eeq where \
$p=\pee{u}{x}$, \ $q=\pee{u}{y}$ \ and \ $u=u(x,y)$. The Pfaffian differential form for \ $u$ \ is \[du=p\:dx+q\:dy.\]
\noi Applying a Legendre like transformation results that
\[d(xp+yq)-du-xdp-ydq=0.\]

\noi As \ $dF=F_pdp+F_qdq=0$, \ then \ \[dp= -(F_q/F_q)dq,\] \
therefore \beq d(xp+yq)-du+\e(x\frac{F_q}{F_p}-y\ri)\:dq=0\:.\eeq

\noi Since this is a Pfaffian differential form then the following theorem can be applied \cite{SNE}:\vs{.4cm}

\textbf{\textit{Theorem}}\vs{.1cm}

\ \ \ \textit{A necessary and sufficient condition that the Pfaffian differential equation \ $\vx\cdot \vr=0$ \ should be integrable is that \ $\vx\cdot rot\vx=0$.}\vs{.3cm}

Therefore from the theorem the condition of integrability applied to \eq (2) results in

\[\vx\cdot rot\vx=-\e(\frac{\pr}{\pr (xp+yq)}+\frac{\pr}{\pr
u}\ri)\e(x\frac{F_q}{F_p}-y\ri)=0,\]

\noi which integrated gives  \beq u-xp-yq=\phi(q).\eeq

\noi The use of this result in \eq (2) supplies the additional condition \beq
\e(x\frac{F_q}{F_p}-y\ri)=-\phi'(q).\eeq

Then the general solution of de PDE is given by the \eq (3) where
\ $q$ \ is determined by \eq  (4) for every choice of the arbitrary \f \ $\phi(q)$. This is a general solution since it has an arbitrary \f  \ $\phi(q)$, i.e., for each form of \ $\phi(q)$ \ the \eqs (4) and (1) results in a system of algebraic \eqs that determines \ $q=q(x,y)$ \ and \
$p=p(x,y)$, which gives a particular solution to the PDE when substituted in \eq (3).

 In some problems it can be writen explicitly \ $p=f(q)$ (or \ $q=g(p)$) \ and the general solution of the PDE
from \eq (3) now stay as  \beq u=x\:f(q)+yq+\phi(q), \eeq

\noi and the integration condition - \eq (4) - turns in \beq xf_q+y=-\phi'(q),\eeq \noi which determines the variable  \ $q=q(x,y)$ \ for every choice of the arbitrary \f \ $\phi(q)$.

Let consider, as a first example, the \eq \ \[F(p,q)= p-Bq+A=0,\] where \ $p$ \ can turn out explicit and \ $A,B$ \ are constants. Therefore \ $p=f(q)=Bq-A$, and \ $f_q=B$ \ then the solution given by \eq (5) is
\[u=xp+yq-\phi(q)=x(Bq-A)+yq-\phi(q),\]\noi where it must be obtained \ $q=q(x,y)$ \ from \eq (6).

The \eq (6) gives \ \[xf_q+y=Bx+y=\phi'(q),\] then \
\[q=\psi(Bx+y)\] \ and the general solution is
\[u=Ax +(Bx+y)\:\psi(Bx+y)-\phi^{-1}(\psi(Bx+y))=Ax+\Phi(Bx+y).\]\noi This is the same result as that given by the method for a Legendre PDE \cite{SNE}.

As another example consider the non linear PDE \ $p^n-q^m=A$,\
where \ $A$ \ is constant. Then \
$p=(A+q^m)^{1/n}=f(q)$, and \
\[f_q=\frac{m}{n}(A+q^m)^{(1-n)/n}q^{m-1}\] therefore from \eq (5) the general solution is \[u=x(A+q^m)^{1/n}+yq-\phi(q),\]
\noi with  \ $q=q(x,y)$ \ obtained from \eq (6) rewritten as
\[\phi'(q)=\frac{mx}{n}(A+q^m)^{(1-n)/n}q^{m-1}+y.\]
\noi Therefore each arbitrary choice of \ $\phi(q)$ \ provides an integral solution.

\section{General Solution to the PDE \ $\mathbf{F(f(x)\:p,q)=0}$ \ (or \ $\mathbf{F(p,h(y)\:q)=0}$)}

\hp Following an identical procedure as that of the last section, from the PDE \ $F(f(x)p,q)=0$ \ it can be written as \ $p=G(q)/f(x)$. The substitution of this result in the differential form associated and the application of a Legendre like transformation results in
\[du=d[H(x)\:G(q)+yq]-[H(x)\:G\:'(q)+y]dq,\]\noi where \ $H(x)=x/f(x)$.

As this is a Pfaffian differential form then the same theorem of the last section can be applied, resulting in the condition of integrability  \beq H(x)G\:'(q)+y=\phi'(q).\eeq

Therefore the general solution is \beq u=H(x)G(q)+yq-\phi(q),\eeq
\noi where \ $\phi(q)$ \ is an arbitrary function, which once selected provides the value of the variable \ $q=q(x,y)$ \ by the \eq (7).

In similar manner it can be obtained the general solution of the PDE \[F(p,h(y)\:q)=0.\] If
\ \[q=G(p)/h(y)\]  then the general solution is given by \[u=xp+G(p)H(y)-\phi(p),\]
\noi where \ $H(y)=y/h(y)$, \ with the integrability condition \[\phi'(p)=G\:'(p)H(y)+x.\]

\section{General Solution To Unidimensional Hamilton-Jacobi Equation}

\hp Considere the most general \hj \eq for a conservative unidimensional non relativistic mechanical system
 \beq a(x)p^2+V(x)-q=0,\eeq \noi where \  $p=\pee{S}{x}$ \ and \ $q=\pee{S}{t}$.

As \ $S=S(x,t)$ then the differential form can be written as \beq dS=pdx+qdt=d(px+qt)-xdp-tdq,\eeq where a Legendre like transformation was applied.

The use of \ $p$ \  from  (9) in the above \eq results in \beq
dS=d\left(\frac{x\sqrt{a(q-V)}}{a}+qt\right)-
\frac{x(a'V-aV'-qa')}{2\sqrt{a(q-V)}}dx-\left( t+\frac{x}{2\sqrt{a(q-V)}}\right) dq,\eeq
where \ $a'=da/dx$ \ and \ $V'=dv/dx$, yielding \beq
S(x,t)=x\sqrt{(q-V)/a}+qt-F(x,q),\eeq sendo \ $F$ \ tal que \beqn
\pe{F}{q}&=&t+\frac{x}{2\sqrt{a(q-V)}}\:,\\
\pe{F}{x}&=&\frac{x(a'V-aV'-qa')}{2\sqrt{a(q-V)}}\,\equiv\,
H(x,q)\:.\eeqn

The integration of the \eq (14) furnishes \ $F=\int H(x,q)dx+G(q)$,  where \
$G$ \ is an arbitrary function. This result applied in \eq (13) gives an \eq that defines the variable \ $q=q(x,t)$, for every arbitrary choice of the function \ $G$: \beq \int \pe{H}{q}dx +
G'(q)=y+\frac{x}{2\sqrt{a(q-V)}}\:.\eeq

As this solution contains an arbitrary \f then \ $S=S(x,t)$ \ given by (12) is a general solution.

Its interesting emphasize that to obtain the solution  by the method of separation of variables to this \hj \eq it is imposed that \ $q=constante$ \ (i.e., $dq=0,\,\,S(x,t)=W(x)+C(t)$).

\section{Examples}

\hp As the first example consider the \hj \eq  that describes a free particle
\ $ap^2-q=0$ \ ($a=constante$). The solution from (12) is \[S=x\sqrt{q/a}+qt-F.\]
\noi where the \f \ $F$ \ is obtained from the solution of the system composed by  (13) e (14)
\[ F\:'(q)=t+\frac{x}{2\sqrt{aq}}\:.\] The last \eq furnishes the variable \
$q=q(x,t)$ \ for every choice of the arbitrary
\f \ $F$. For example, if \ $F=Cq$ \ then \ $q=x^2/4a(C-t)^2$, \
therefore \ $S(x,t)=x^2/4a(C-t)$. \ This is the same solution obtained using the movement data of the
particle \cite{SAL}, which is unnecessary in our method.

The solution \ $S(x,t)=x\sqrt{C/a}+Ct$ \ given by method of separation of variables applied to this
\hj \eq is obtained making \ $dq=0$ \ in (11).

Let consider as another example the \hj \eq of a simple harmonic oscillator
\[p^2+x^2-q=0.\]

The solution from (12) is \[S(x,t)=qt+\frac{x}{2}\sqrt{q-x^2}+\frac{q}{2}\sin^{-1} \frac{x}{\sqrt{q}}-G(q),\]
 \noi where \ $q=q(x,t)$ \ is fixed by each choice of the arbitrary \f  \ $G$ \ from
\[t+\frac{1}{2}x\sqrt{q-x^2}=G'(q)+\frac{q}{2}\:sen^{-1}\left(\frac{x}{\sqrt{q}}\right).\]
The solution form the method of separation of variables ($q=C$) \ is
\[S(x,t)=\frac{1}{2}x\sqrt{C-x^2}+Ct-\frac{C}{2}\:sen^{-1}\left(\frac{x}{\sqrt{C}}\right).\]

\section {FINAL REMARKS}

\hp The \pro developed to solve the one dimensional \hj \eq is an extension of that presented in section 2 and 3 \cite{ENA,CNM}. The integrability condition for Pfaffian differential forms (3) \cite{SNE} imply in the \eqs
(13) and (14).

The  extension of this method to \hj \eqs two and tridimensional and a general formulation for this type of PDEs can be a later approach.

\vs{.8cm}
\noi {\bf Acknowledgments} \vs{.2cm}

The author is grateful to Dr. Oslim Espindola (in memoriam) and to Dr. Nelson Lima Teixeira (in
memoriam) profitable debates.\vs{.8cm}

\end{document}